\documentclass[12pt]{article}

\usepackage{amssymb}
\usepackage{amsmath}
\usepackage{cite}

\begin{document}

\renewcommand{\citeleft}{{\rm [}}
\renewcommand{\citeright}{{\rm ]}}
\renewcommand{\citepunct}{{\rm,\ }}
\renewcommand{\citemid}{{\rm,\ }}

\newcounter{abschnitt}
\newtheorem{satz}{Theorem}
\newtheorem{theorem}{Theorem}[abschnitt]
\newtheorem{koro}[satz]{Corollary}
\newtheorem{prop}[theorem]{Proposition}
\newtheorem{lem}[theorem]{Lemma}
\newtheorem{definition}[theorem]{Definition}
\newtheorem{conj}[theorem]{Conjecture}
\newtheorem{corol}[theorem]{Corollary}

\newcounter{saveeqn}
\newcommand{\alpheqn}{\setcounter{saveeqn}{\value{abschnitt}}
\renewcommand{\theequation}{\mbox{\arabic{saveeqn}.\arabic{equation}}}}
\newcommand{\reseteqn}{\setcounter{equation}{0}
\renewcommand{\theequation}{\arabic{equation}}}

\hyphenation{convex} \hyphenation{bodies}

\sloppy

\phantom{a}

\vspace{-1.7cm}

\begin{center}
 \begin{large} {\bf Volume Inequalities for Asymmetric Wulff Shapes} \\[0.6cm] \end{large}

\begin{large} Franz E. Schuster and Manuel Weberndorfer \end{large}
\end{center}

\vspace{-1cm}

\begin{quote}
\footnotesize{ \vskip 1truecm\noindent {\bf Abstract.} Sharp
reverse affine isoperimetric inequalities for asymmetric Wulff
shapes and their polars are established, along with the
characterization of all extremals. These new inequalities have as
special cases previously obtained simplex inequalities by Ball,
Barthe and Lutwak, Yang, and Zhang. In particular, they provide
the solution to a problem by Zhang.}
\end{quote}

\vspace{0.7cm}

\centerline{\large{\bf{ \setcounter{abschnitt}{1}
\arabic{abschnitt}. Introduction}}} \alpheqn

\vspace{0.7cm}

Over the last decades considerable progress has been made in
establishing reverse (affine) isoperimetric inequalities, that is,
inequalities which usually have simplices or, in the symmetric
case, cubes and their polars, as extremals. By the end of the
1980s, only a very small number of significant reverse
inequalities had been obtained and no systematic approach towards
these inequalities seemed within reach. A breakthrough occurred
when Ball \textbf{\cite{ball89, ball91b}} \linebreak discovered a
reformulation of the Brascamp--Lieb inequality by exploiting the
notion of {\it isotropic measure} which, in turn, is connected to
a variety of \linebreak extremal problems in geometric analysis
(see \textbf{\cite{giannomil, gruber08, J-extproineassubcon,
LYZ2005, petty61}}). Ball's geometric Brascamp--Lieb inequality
was tailor-made to establish several \linebreak important new
reverse inequalities.

Settling the uniqueness of the extremals for the newly obtained
reverse isoperimetric inequalities with an underlying {\it
discrete} isotropic measure was  made possible only through an
optimal transport approach of Barthe \textbf{\cite{barthe97,
barthe98i}} \linebreak towards establishing not only the
Brascamp--Lieb inequality but also its \linebreak inverse form
conjectured by Ball. In this way, for example, sharp $L_p$ volume
ratio inequalities and their duals were established (see Section
6 for details and further examples). To obtain uniqueness of
extremals when the isotropic measure underlying the extremal
problem is not necessarily discrete, Lutwak, \linebreak Yang, and
Zhang \textbf{\cite{LYZ2004, LYZ2007, LYZ2009}} developed a new
approach based on the \linebreak Ball--Barthe techniques -- but
not on the Brascamp--Lieb inequality or its inverse --
 which allowed them to extend all the new reverse inequalities
to the setting of general isotropic measures along with
characterizations of all extremizers. Later, Barthe
\textbf{\cite{barthe04}} established a continuous version of the
Brascamp--Lieb inequality and its inverse which also yields these
general \linebreak reverse inequalities along with their equality
conditions.

The notion of Wulff shapes has its origins in the classical theory
of crystal growth. In more modern mathematical terms it provides
a unifying setting \linebreak for several extremal problems with
an underlying isotropic measure. Sharp reverse volume inequalities
for  {\it origin-symmetric} Wulff shapes and their polars were
obtained by Lutwak, Yang, and Zhang \textbf{\cite{LYZ2004}} (see
Section 3). These inequalities generalize several of the
previously obtained Ball--Barthe volume ratio inequalities for
unit balls of subspaces of $L _p$. The problem of finding similar
volume estimates for not necessarily origin-symmetric Wulff shapes
remained open.

In this article we establish such sharp reverse isoperimetric
inequalities for asymmetric Wulff shapes and their polars
including a complete description of all equality cases. Special
cases of these new Wulff shape inequalities are previously
obtained simplex inequalities of Ball \textbf{\cite{ball91b}},
Barthe \textbf{\cite{barthe97}} and Lutwak, Yang, and Zhang
\textbf{\cite{LYZ2009, LYZ2007, LYZ-newellassconbod}}.

\vspace{0.2cm}

The setting for this article is Euclidean space $\mathbb{R}^n$,
$n \geq 2$. A convex body is a compact convex set and in this
article will always be assumed to contain the origin in its
interior. The polar body of a convex body $K$ is given by
$K^*=\{x \in \mathbb{R}^n: x \cdot y \leq 1 \mbox{ for all } y
\in K\}$. Throughout, all Borel measures are understood to be
non-negative and finite. We write $\mathrm{supp}\,\nu$ for the
support of a measure $\nu$ and we use $\mathrm{conv}\,L$ to
denote the convex hull of a set $L$.

The main objects of this paper are Wulff shapes. This notion was
introduced at the turn of the previous century by Wulff
\textbf{\cite{W-zurfradergesdeswacundderaufderkry}}, who
conjectured that this shape describes the minimizer of the
interfacial free energy among a crystal's possible shapes of
given volume (see also, e.g., \textbf{\cite{S-thecrygromor}}).
Variations of Wulff's original definition (expanding the class of
admitted parameters) yield versatile geometric objects that have
been analyzed extensively (see e.g. \linebreak
\textbf{\cite{figamagg10, gardner02, HLYZ-eveorlminpro}} or
\textbf{\cite[\rm Section 6.5]{schneider93}}).

\vspace{0.3cm}

\noindent \textbf{Definition} \emph{Suppose $\nu$ is a Borel
measure on $S^{n-1}$ and $f$ is a positive \linebreak continuous
function on $S^{n-1}$. The \emph{Wulff shape} $W_{\nu,f}$
determined by $\nu$ and $f$ is defined by}
\[W_{\nu,f} :=\left\{x \in \mathbb{R}^n: x \cdot u \leq f(u) \text{ for all } u \in \mathrm{supp}\,\nu \right\}.\]

\vspace{0.1cm}

Without further assumptions on the measure $\nu$ and the function
$f$, Wulff shapes, while always convex, may be unbounded. In
order to guarantee that $W_{\nu,f}$ is a convex body, we
introduce the notion of \emph{$f$-centered} measures and consider
only Wulff shapes determined by measures $\nu$ which are
\emph{$f$-centered} and \emph{isotropic}.

\pagebreak

Let $f$ be a positive continuous function on $S^{n-1}$. A Borel
measure $\nu$ on $S^{n-1}$ is called \emph{$f$-centered} if
\[\int_{S^{n-1}} f(u)\,u\, d\nu(u)=o.\]
The measure $\nu$ is called \emph{isotropic} if
\begin{equation} \label{defisotr}
\int_{S^{n-1}} u \otimes u \,d\nu(u)=I_n,
\end{equation}
where $u \otimes u$ is the orthogonal projection onto the line
spanned by $u$ and $I_n$ denotes the identity map on
$\mathbb{R}^n$.

In order to establish reverse affine isoperimetric inequalities,
it is often \linebreak critical to exploit special positions of
convex bodies which, in turn, are \linebreak characterized by
isotropic measures (see, e.g., \textbf{\cite{ball91b, barthe97,
barthe98, barthe98i, giannomil, LYZ-newellassconbod, LYZ2005,
LYZ2009, mareschus}}). The notion of $f$-centered isotropic
measures is designed to unify approaches towards reverse
inequalities which are based on the isotropic-embedding
techniq\-ue introduced by Lutwak, Yang, Zhang
\textbf{\cite{LYZ2009}} (see Section 4 for detail\-s).

Sharp volume estimates for {\it origin-symmetric} Wulff shapes
and their \linebreak polars determined by even functions and
isotropic measures (clearly, they are $f$-centered) can easily be
deduced from previous work of Lutwak, Yang, and Zhang
\textbf{\cite{LYZ2004}}, see Section 3. The extremal
configurations here are given by constant functions and measures
for which the convex hull of their support is a cube. The natural
problem to determine sharp volume bounds for not necessarily
symmetric Wulff shapes was posed by Zhang \textbf{\cite{zhang05}}.

With our first main result we establish such a sharp bound for
the \linebreak volume of the Wulff shape $W_{\nu,f}$ determined
by an $f$-centered isotropic \linebreak measure $\nu$. It depends
on the {\em displacement} of $W_{\nu,f}$ defined by
\[\mathrm{disp}\, W_{\nu,f} := \mathrm{cd}\, W_{\nu,f} \cdot \int_{S^{n-1}} \frac{u}{f(u)}\, d\nu(u),\]
where $\mathrm{cd}\, W_{\nu,f}$ denotes the centroid of the
convex body $W_{\nu,f}$.

\vspace{0.2cm}

\begin{satz} \label{main1}
  Suppose $f$ is a positive continuous function on $S^{n-1}$ and $\nu$ is an isotropic $f$-centered measure.
  If $\mathrm{disp}\, W_{\nu,f}=0$, then
\begin{equation*}
V(W_{\nu,f}) \leq \frac{ (n+1)^{(n+1)/2}}{n!}\| f\|^{n}_{L^2(\nu)}
\end{equation*}
with equality if and only if $\mathrm{conv}\,\mathrm{supp}\,\nu$
is a regular simplex inscribed in $S^{n-1}$ and $f$ is constant
on $\mathrm{supp}\,\nu$.
\end{satz}

\vspace{0.2cm}

In fact, the assumption $\mathrm{disp}\, W_{\nu,f}=0$ is not
necessary. Our proof of Theorem 1 yields an explicit description
of how the displacement enters the sharp upper bound for the
volume of $W_{\nu,f}$ (see Theorem 5.1).

Our second main result, a natural dual to Theorem \ref{main1},
provides a sharp lower bound for the volume of the polar of the
Wulff shape $W_{\nu,f}$. Note that it is independent of the
displacement of the Wulff shape.

\vspace{0.1cm}

\begin{satz} \label{main2}
Suppose $f$ is a positive continuous function on $S^{n-1}$ and
$\nu$ is an isotropic $f$-centered measure. Then
\[V\left(W_{\nu,f}^* \right) \geq \frac{(n+1)^{(n+1)/2}}{n!} \|f\|_{L^2(\nu)}^{-n}\]
with equality if and only if $\mathrm{conv}\,\mathrm{supp}\,\nu$
is a regular simplex inscribed in $S^{n-1}$ and $f$ is constant on
$\mathrm{supp}\,\nu$.
\end{satz}

\vspace{0.1cm}

Our proofs of Theorems 1 and 2 are based on a refinement of the
\linebreak approach towards recently established simplex
inequalities by Lutwak, Yang, and Zhang \textbf{\cite{LYZ2007,
LYZ2009}}, which, in turn, uses many ideas of Ball and Barthe. We
remark, however, that our results can also be obtained by
applications of Barthe's continuous Brascamp--Lieb inequality and
its inverse. The more direct approach we have chosen has the
advantage to be at the same time reasonably self contained and
elementary.

In Section 6 we show how Theorems 1 and 2 directly imply a number
of reverse isoperimetric inequalities (including all equality
conditions) obtained by Ball \textbf{\cite{ball91b, ball89}},
Barthe \textbf{\cite{barthe97}} and Lutwak, Yang, Zhang
\textbf{\cite{LYZ2004, LYZ2007, LYZ2009}}. As an example, we
state here one corollary of Theorem \ref{main2}. It was first
established in \textbf{\cite{LYZ2009}} and provides a lower bound
for the volume of the polar of a convex body $K$ in terms of the
volume of its dual Legendre ellipsoid $\Gamma_{-2}K$ introduced
in \textbf{\cite{LYZ-newellassconbod}} (see Section 2 for precise
definitions).

\vspace{0.1cm}

\begin{koro} \label{introgam2} If $K$ is a convex body in $\mathbb{R}^n$, then
\[V(K^*)V(\Gamma_{-2}K) \geq  \frac{\kappa_n (n+1)^{(n+1)/2}}{n!n^{n/2}}   \]
with equality if and only if $K$ is a simplex whose centroid is
at the origin.
\end{koro}

\vspace{0.1cm}

Here and in the following, $\kappa_n$ denotes the volume of the
Euclidean unit ball in $\mathbb{R}^n$.

\pagebreak

\centerline{\large{\bf{ \setcounter{abschnitt}{2}
\arabic{abschnitt}. Background material}}} \reseteqn \alpheqn

\vspace{0.6cm}

For quick later reference, we collect in this section the
necessary background material. In particular, we list basic
auxiliary facts from the $L_p$ Brunn--Minkowski theory and recall
a number of special positions of convex bodies (needed in the last
section). As a general reference, the reader may wish to consult
the books \textbf{\cite{gardner95, schneider93}} and the articles
\textbf{\cite{L-brutheimixvolminpro, L-brutheiiaffgeosurare}}.

Throughout we will denote by $e_1, \ldots, e_n$ the standard
Euclidean basis of $\mathbb{R}^n$ and we use $\|\cdot \|$ to
denote the standard Euclidean norm on $\mathbb{R}^n$. We emphasize
that in this article a convex body in $\mathbb{R}^n$ is a compact
convex set that contains the origin in its interior. A convex
body $K$ is uniquely determined by its (positive, sublinear)
support function defined by
\[h(K,x):=\max\{x \cdot y: y \in K\},  \qquad x \in \mathbb{R}^n.\]

A convex body $K$ in $\mathbb{R}^n$ is also determined up to
translation by its surface area measure $S(K,\cdot)$. Recall that
for a Borel set $\omega \subseteq S^{n-1}$, $S(K,\omega)$ is the
$(n-1)$-dimensional Hausdorff measure of the set of all boundary
points of $K$ at which there exists a normal vector of $K$
belonging to $\omega$. It is well known that the surface area
measure of a convex body $K$ is $1$-centered, that is,
\begin{equation}\label{eq:SurfCentred}
\int_{S^{n-1}} u\, d S(K,u)=o.
\end{equation}

If $K$ and $L$ are convex bodies and $\alpha, \beta \geq 0$ (not
both zero), then their $L_p$ Minkowski combination $\alpha \cdot
K \, +_{p} \, \beta \cdot L$ is the convex body whose support
function is given by
\[h(\alpha \cdot K \, +_{p} \, \beta \cdot L,\cdot)^p = \alpha h(K,\cdot)^p + \beta h(L,\cdot)^p.\]
In \textbf{\cite{L-brutheimixvolminpro, L-brutheiiaffgeosurare}},
Lutwak showed that merging the notion of volume with these $L_p$
Minkowski combinations of convex bodies, introduced by Firey,
leads to a Brunn--Minkowski theory for each $p \geq 1$. In
particular, the $L_p$ mixed volume $V_p(K,L)$ was defined in
\textbf{\cite{L-brutheimixvolminpro}} by
\[\frac{n}{p}V_p(K,L)=\lim \limits_{\varepsilon \rightarrow 0^+} \frac{V(K +_p \varepsilon \cdot L)-V(K)}{\varepsilon}.  \]
Clearly, for $K = L$, we have
\begin{equation} \label{vol1}
V_p(K,K)=V(K).
\end{equation}
It was shown in \textbf{\cite{L-brutheimixvolminpro}} that
corresponding to each convex body $K$, there exists a positive
Borel measure on $S^{n-1}$, the $L_p$ {\em surface area measure}
$S_p(K,\cdot)$ of $K$, such that for every convex body $L$,
\begin{equation} \label{defvp}
V_p(K,L)=\frac{1}{n} \int_{S^{n-1}} h(L,u)^p dS_p(K,u).
\end{equation}
The measure $S_1(K,\cdot)$ is just the surface area measure of
$K$.

The $L_p$ surface area measure is absolutely continuous with
respect to $S(K,\cdot)$, more precisely,
\begin{equation}\label{absc}
dS_p(K,u)=h(K,u)^{1-p}\,dS(K,u), \qquad u \in S^{n-1}.
\end{equation}

From (\ref{absc}) and the definition of surface area measures, it
follows easily that, for a given convex body $K$, the Wulff shape
determined by the $L_p$ surface area measure $S_p(K,\cdot)$ and
the support function $h(K,\cdot)$ of $K$ is precisely the body
$K$, i.e.,
\begin{equation} \label{eq:KIntersectionHalfSp}
W_{S_p(K,\cdot),h(K,\cdot)}=K.
\end{equation}

A $\mathrm{GL}(n)$ image of a convex body is often called a {\it
position} of the body. Special positions have been the focus of
intensive investigations, in particular, in relation with a
variety of extremal problems for geometric invariants of the
bodies in special position (see e.g. \textbf{\cite{ball91b,
barthe98, giannomil, gruber08, LYZ2005}}).

A classical example of an important special position of a convex
body $K$ is the {\it John position}: Let $\mathrm{J}K$ denote the
unique ellipsoid of maximal volume contained in $K$. The body $K$
is said to be in John position, if $\mathrm{J}K$ coincides with
the Euclidean unit ball $B$. The following well known
characterization of this position goes back to John
\textbf{\cite{J-extproineassubcon}}:

\begin{prop} \label{johnth} A convex body $K$ which contains the
unit ball $B$ is in John position if and only if there exists an
1-centered isotropic measure on $S^{n-1}$ supported by
contact points of $K$ and $B$.
\end{prop}

A natural dual to the John position of a convex body $K$ is the
{\it Loewner position}, here the ellipsoid of minimal volume
containing $K$ is the unit ball. It was also characterized in
\textbf{\cite{J-extproineassubcon}} by the existence of a
$1$-centered isotropic measure supported by the contact points of
$K$ and $B$.

Another classical position is closely related to the problem of
finding a reverse form of the Euclidean isoperimetric inequality.
Since convex bodies of a given volume may have arbitrarily large
surface area, it is natural to consider convex bodies in {\it
minimal surface area position}, that is, the surface area of the
bodies is minimal among all their affine images of the same
volume. Petty \textbf{\cite{petty61}} showed that a convex body
$K$ is in minimal surface area position if and only if its surface
area measure $S(K,\cdot)$ is isotropic (up to scaling).

In a more recent article, Lutwak, Yang, Zhang
\textbf{\cite{LYZ2005}} have shown that the John position and the
minimal surface area position are in fact special cases ($p =
\infty$ and $p = 1$) of a family of $L_p$ John positions of a
given convex body.

\vspace{0.3cm}

\noindent \textbf{Definition (\!\!\cite{LYZ2005})} \emph{Suppose
$K$ is a convex body and $1 \leq p \leq \infty$. Amongst all
origin-symmetric ellipsoids $E$, the unique ellipsoid that solves
the constrained extremal problem
\[\max_EV(E) \quad \mbox{subject to} \quad V_p(K,E) \leq V(K)\]
will be called the \emph{$L_p$ John ellipsoid} $\mathrm{E}_pK$ of
$K$. We say $K$ is in \emph{$L_p$ John position} if
$\mathrm{E}_pK$ coincides with the Euclidean unit ball $B$.}

\vspace{0.3cm}

The $L_1$ John ellipsoid of a convex body $K$ is also called the
Petty ellipsoid. It is not difficult to show (cf.\
\textbf{\cite{LYZ2005}}) that (up to scaling) $K$ is in $L_1$
John position if and only if $K$ is in minimal surface area
position. The $L_{\infty}$ John ellipsoid is the origin-centered
ellipsoid of maximal volume contained in $K$. Hence,
$\mathrm{E}_\infty K=\mathrm{J}K$ if the John ellipsoid of $K$ is
centered at the origin.

Of particular importance among the family of $L_p$ John
ellipsoids of a given body $K$ is also the $L_2$ John ellipsoid.
This ellipsoid was previously discovered by Lutwak, Yang, and
Zhang (see \textbf{\cite{LYZ-newellassconbod, LYZ-crainestabod}})
and denoted by $\Gamma_{-2}K$. This notation should indicate a
duality with the classical Legendre ellipsoid $\Gamma_2K$. In
fact Ludwig \textbf{\cite{L-ellmatval}} showed that the operators
$\Gamma_{-2}$ and $\Gamma_{2}$ are the only linearly intertwining
maps on convex bodies that satisfy the inclusion-exclusion
principle (see also \textbf{\cite{haberljems, Ludwig:Minkowski,
L-minareval, ludwig11}} for related results).

The following characterization of the $L_p$ John position of a
convex body in terms of isotropic measures was also established
in \textbf{\cite{LYZ2005}}:

\begin{prop} \label{lpjohnth} Suppose that $p \geq 1$. A convex body $K$
is in $L_p$ John position if and only if its $L_p$ surface area
measure $S_p(K,\cdot)$ is isotropic up to volume normalization.
\end{prop}

We conclude this section with another auxiliary result
\textbf{\cite[\textnormal{Theorem 5.1}]{LYZ2005}} concerning
monotonicity properties of $L_p$ John ellipsoids: If $K$ is a
convex body and $0<p\leq q\leq \infty$, then
\begin{equation}\label{eq:VIneqpq}
    V(\mathrm{E}_qK)\leq V(\mathrm{E}_pK).
\end{equation}

\pagebreak

\centerline{\large{\bf{ \setcounter{abschnitt}{3}
\arabic{abschnitt}. Volume inequalities for symmetric Wulff shapes}}}
\reseteqn \alpheqn \setcounter{theorem}{0}

\vspace{0.6cm}

In this section we state the volume inequalities corresponding to
our main results when the considered Wulff shapes are all
origin-symmetric, that is, they are determined by even isotropic
measures and even functions. In a slightly different formulation
these inequalities were established by Lutwak, Yang, and Zhang
\textbf{\cite{LYZ2004}} and generalize previous results by Ball
and Barthe. We also sketch a proof of one of these inequalities
using Ball's geometric \linebreak Brascamp--Lieb inequality, in
order to emphasize the close connection \linebreak between this
analytic inequality and volume inequalities for (symmetric as
well as asymmetric) Wulff shapes.

\vspace{0.4cm}

\noindent {\bf The Brascamp--Lieb Inequality.}
\emph{Suppose that $u_1,\ldots, u_m \in S^{n-1}$ and $c_1,
\ldots, c_m
> 0$ such that
\[\sum_{i=1}^m c_i u_i \otimes u_i=I_n.   \]
If  $g_i:\mathbb{R} \to [0,\infty)$, $1 \leq i \leq m$, are
integrable functions, then}
\[\int_{\mathbb R^n} \prod_{i=1}^m g_i(u_i\cdot x)^{c_i} \, dx \leq \prod_{i=1}^m \left(\int_{\mathbb R} g_i\right)^{c_i}.  \]

\vspace{0.4cm}

The Brascamp--Lieb inequality \textbf{\cite{braslieb}} was
established to prove the sharp form of Young's convolution
inequality. Around 1990 Ball \textbf{\cite{ball89}} discovered
the \linebreak geometric reformulation stated above (later
generalized by Barthe \textbf{\cite{barthe98i, barthe04}}) which
allowed a simple computation of the optimal constant. It directly
yields the following sharp volume bound for $W_{\nu,f}$, when the
underlying isotropic measure is even and discrete.

\begin{theorem} \label{th:evenWulff} {\bf (\!\!\cite{LYZ2004})}
Suppose $f$ is an even positive continuous function on $S^{n-1}$ and $\nu$ is an even isotropic measure.
Then
\begin{equation} \label{gandalf}
V(W_{\nu,f}) \leq \left(\frac{2}{\sqrt n}\right)^n \|f\|_{L^2(\nu)}^n
\end{equation}
with equality if and only if $\mathrm{conv}\,\mathrm{supp}\,\nu$ is a cube inscribed in $S^{n-1}$ and $f$ is constant on $\mathrm{supp}\,\nu$.
\end{theorem}

\pagebreak

\noindent {\it Sketch of the proof.} We prove inequality
(\ref{gandalf}) for the case of a discrete measure $\nu$
supported, say, on $\pm\, u_1,\dots,\pm\, u_m \in S^{n-1}$. The
Wulff shape $W_{\nu,f}$ is thus given by
\begin{equation} \label{symproof1}
W_{\nu,f}=\bigcap_{i=1}^m\left\{x\in\mathbb R^n:\;|x\cdot u_i|\leq f(u_i) \right\}.
\end{equation}
For $1 \leq i \leq m$, let $c_i > 0$ be the total mass of $\nu$
at the points $\pm\,u_i$ and define the function $g_i: \mathbb R
\rightarrow [0,\infty)$ by
\begin{equation} \label{symproof2}
g_i(t)= \mathbb{I}_{[-f(u_i),f(u_i)]}(t).
\end{equation}
By (\ref{symproof1}), (\ref{symproof2}) and the Brascamp--Lieb inequality, we now obtain
\[V(W_{\nu,f})= \int_{\mathbb R^n}\mathbb{I}_{W_{\nu,f}}(x)\,dx
         =\int_{\mathbb R^n} \prod_{i=1}^m g_i(x\cdot u_i)^{c_i} \, dx
         \leq 2^{\sum_{i=1}^mc_i} \prod_{i=1}^m f(u_i)^{c_i}.\]
Since the measure $\nu$ is isotropic, taking traces in
(\ref{defisotr}) shows that $\sum_{i=1}^mc_i=n$. Consequently, an
application of the arithmetic-geometric mean inequality yields
\[V(W_{\nu,f}) \leq 2^n \left(\frac{1}{n}\sum_{i=1}^m c_if(u_i)^2\right)^{n/2}
            = \left(\frac{2}{\sqrt n}\right)^n \|f\|_{L^2(\nu)}^n.  \]

\vspace{-0.3cm}

\hfill $\blacksquare$

\vspace{0.3cm}

The following result is dual to Theorem \ref{th:evenWulff}; for
the case of even and discrete isotropic measures, it follows from
Barthe's inverse Brascamp--Lieb inequality
\textbf{\cite{barthe98i}} (by arguments similar to the ones
sketched above).

\begin{theorem} \label{th:evenpolWulff} {\bf (\!\!\cite{LYZ2004})}
 Suppose $f$ is an even positive continuous function on $S^{n-1}$ and $\nu$ is an even isotropic measure.
  Then
 \begin{align*}
  V( W_{\nu,f}^*) &\geq \frac{\left(2 \sqrt n\right)^n}{n!} \|f\|_{L^2(\nu)}^{-n}
 \end{align*}
with equality if and only if $\mathrm{conv}\,\mathrm{supp}\,\nu$ is a cube inscribed in $S^{n-1}$ and $f$ is constant on $\mathrm{supp}\,\nu$.
\end{theorem}

The case $f\equiv 1$ of Theorem \ref{th:evenWulff} was proved by
Ball \textbf{\cite{ball91b}}, the equality conditions for discrete
measures were obtained by Barthe \textbf{\cite{barthe98i}}.
Theorem \ref{th:evenpolWulff} for $f \equiv 1$ and discrete
measures was proved by Barthe \textbf{\cite{barthe98i}}.

\pagebreak

In order to establish Theorems \ref{th:evenWulff} and
\ref{th:evenpolWulff} for general even isotropic \linebreak
measures, Lutwak, Yang, and Zhang \textbf{\cite{LYZ2004}} used a
direct approach based on optimal mass transport and a determinant
inequality, called the Ball--Barthe Lemma (see the next section),
that has easily stated equality cases obtained in
\textbf{\cite{LYZ2004}}. Another possibility towards proving
Theorems \ref{th:evenWulff} and \ref{th:evenpolWulff} is to employ
Barthe's continuous versions of the Brascamp--Lieb inequality and
its inverse, the equality conditions of which are also based on
the Ball--Barthe Lemma. It is therefore no surprise that this
basic inequality is also critical in the proofs of our main
results. Moreover, to demonstrate the extremal property of the
regular simplex in our inequalities, we also need an important
embedding of $f$-centered isotropic measures introduced by
Lutwak, Yang, and Zhang \textbf{\cite{LYZ2007, LYZ2009}} which we
review in the next section.

\vspace{0.8cm}

\centerline{\large{\bf{ \setcounter{abschnitt}{4}
\arabic{abschnitt}. Isotropic embeddings and the Ball--Barthe
Lemma}}} \reseteqn \alpheqn \setcounter{theorem}{0}

\vspace{0.5cm}

In the following we recall the concept of isotropic embeddings
which was introduced by Lutwak, Yang, Zhang
\textbf{\cite{LYZ2009}}. These embeddings lift $f$-centered
isotropic measures on $S^{n-1}$ to isotropic measures on $S^n$ and
at the same time map the vertices of the regular $n$-simplex
inscribed in $S^{n-1}$ to an orthonormal basis in
$\mathbb{R}^{n+1}=\mathbb{R}^n \times \mathbb{R}$. This latter
property ensures that we can apply the Ball--Barthe Lemma to
obtain a sharp bounds in our main results.

\vspace{0.3cm}

\noindent \textbf{Definition} \emph{If $\nu$ is a Borel measure on
$S^{n-1}$, then a continuous function $g: S^{n-1} \to \mathbb
R^{n+1} \setminus \{o\}$ is called an \emph{isotropic embedding}
of $\nu$ if the measure $\overline\nu$ on $S^n$, defined by
\begin{equation} \label{isoemb}
\int_{S^{n}} t(w)\, d \overline\nu(w)=\int_{S^{n-1}}
t\left(\frac{g(u)}{\|g(u)\|}\right)\|g(u)\|^2d\nu(u)
\end{equation}
for every continuous $t: S^n \rightarrow \mathbb{R}$, is
isotropic.}

\vspace{0.3cm}

Of particular interest for us are isotropic embeddings of already
isotropic measures. A natural class of such embeddings can be
characterized in the following way.

\begin{lem}\label{le:IsotrChar}
Suppose $f$ is a positive continuous function on $S^{n-1}$ and
$\nu$ is an isotropic measure on $S^{n-1}$. Then $g_\pm: S^{n-1}
\rightarrow \mathbb{R}^{n+1}=\mathbb{R}^n \times \mathbb{R}$,
defined by
\begin{equation} \label{impisoemb}
g_\pm(u)=(\pm u,f(u)),
\end{equation}
are isotropic embeddings of $\nu$ if and only if $\nu$ is
$f$-centered and $\|f\|_{L_2(\nu)}=1$.
\end{lem}

\pagebreak

\noindent {\it Proof.} If $\overline\nu$ is defined as in
(\ref{isoemb}), where $g$ is replaced by $g_\pm$, we have
\[\int_{S^{n}} w\otimes w\, d \overline\nu(w)=\int_{S^{n-1}}
\begin{pmatrix} u\otimes u & \pm f(u)u \\ \pm f(u)u^{\mathrm{T}} & f^2(u) \end{pmatrix}\, d
\nu(u). \]
   Consequently, since $\nu$ is isotropic,
   \[\int_{S^{n}} w\otimes w\, d \overline\nu(w)= I_{n+1} \]
if and only if $\nu$ is $f$-centered and $\|f\|_{L_2(\nu)} = 1$.
\hfill $\blacksquare$

\vspace{0.4cm}

Note that for isotropic embeddings of the form (\ref{impisoemb})
the last coordinate (with respect to the decomposition
$\mathbb{R}^{n+1} = \mathbb{R}^n \times \mathbb{R}$) of all the
points in the support of the measure $\overline{\nu}$, defined by
(\ref{isoemb}), is positive.

The following two special cases of isotropic embeddings of the
form (\ref{impisoemb}) have played a critical role in the proof
of a number of reverse isoperimetric inequalities having
simplices as extremals (see \textbf{\cite{ball91b, barthe98,
LL-meawidineisomea, LYZ2007, LYZ2009}}).

\vspace{0.3cm}

\noindent {\bf Examples:}
\begin{enumerate}
\item[(a)] If $\nu$ is an 1-centered isotropic measure (e.g., the normalized surface area measure of a convex body in minimal surface area position), then $g_\pm:
S^{n-1} \rightarrow \mathbb{R}^{n+1}$, defined by
\begin{equation} \label{constemb}
g_{\pm}(u)=\left(\pm u,\frac{1}{\sqrt n}\right),
\end{equation}
are isotropic embeddings of $\nu$.

\item[(b)] Suppose that $K$ is a convex body in $L_2$ John
position. Since, by (\ref{eq:SurfCentred}) and (\ref{absc}),
\[\int_{S^{n-1}} h(K,u) \, d  S_2(K,u)=o\]
and, by (\ref{vol1}) and (\ref{defvp}),
\[\frac{1}{V(K)}\int_{S^{n-1}} h(K,u)^2\, dS_2(K,u)=n,\]
it follows from Proposition \ref{lpjohnth} that $g_\pm: S^{n-1}
\rightarrow \mathbb{R}^{n+1}$, defined by
\[g_{\pm}(u)=\left(\pm u, \frac{h(K,\cdot)}{\sqrt n} \right)\]
are isotropic embeddings of $S_2(K,\cdot)/V(K)$.
\end{enumerate}

\pagebreak

The approaches towards sharp reverse isoperimetric inequalities of
both Ball, Barthe and Lutwak, Yang, Zhang make critical use of the
following basic estimate for the determinant of a weighted sum of
rank-one projections:

\vspace{0.4cm}

\noindent {\bf The Ball--Barthe Lemma.} \emph{If $\overline{\nu}$
is an isotropic measure on $S^n$ and $t$ is a positive continuous
function on $\mathrm{supp}\,\overline{\nu}$, then
\begin{equation} \label{bblem}
\det \int_{S^n} t(w)w\otimes w\, d \overline{\nu}(w)\geq
\exp\left( \int_{S^n} \log t(w) \, d\overline{\nu}(w)\right),
\end{equation}
with equality if and only if $t(v_1)\cdots t(v_{n+1})$ is constant
for linearly independent $v_1, \ldots, v_{n+1} \in
\mathrm{supp}\,\overline{\nu}$.}

\vspace{0.4cm}

For discrete measures, inequality (\ref{bblem}) goes back to
Ball. In \textbf{\cite{barthe98i}} Barthe provides a simple
proof. The equality conditions for (\ref{bblem}) were obtained
using mixed discriminants and H\"older's inequality by Lutwak,
Yang, Zhang \textbf{\cite{LYZ2004}}.

The Ball--Barthe Lemma also plays a crucial role in the proof of
our main \linebreak results, in particular, for establishing the
equality cases. Our next lemma goes back to arguments employed by
Lutwak, Yang, Zhang \textbf{\cite{LYZ2007}}. It uses the equality
conditions for (\ref{bblem}) to characterize the support of
1-centered isotropic measures which are embedded by the functions
given in (\ref{constemb}).

\begin{lem}\label{equcond}
Let $\nu$ be an $1$-centered isotropic measure on $S^{n-1}$, let
$\overline{\nu}_{\pm}$ denote the isotropic measures on $S^n$
defined by (\ref{isoemb}), isotropically embedded by $g_{\pm}$
defined in (\ref{constemb}), and let $D \subseteq
\mathbb{R}^{n+1}$ be an open cone with apex at the origin
containing $e_{n+1}$ such that $w \cdot z
> 0$ for every $w \in \mathrm{supp}\,\overline\nu$ and $z \in D$.

For every $z \in D$, define $t_{z}: \mathrm{supp}\,\overline{\nu}
\rightarrow (0,\infty)$ by
\[t_z(w)=\phi_w(w\cdot z),  \]
where $\phi_w: (0,\infty) \rightarrow (0,\infty)$ is smooth
nonconstant and depends continuously on $w$. If there is equality
in (\ref{bblem}) for $\overline{\nu}_+$, or $\overline{\nu}_-$
respectively, and every $t_{z}$, $z\in D$, then $\mathrm{conv}\,
\mathrm{supp}\, \nu$ is a regular simplex inscribed in $S^{n-1}$.
\end{lem}

\noindent {\it Proof.} We prove the statement for
$\overline{\nu}_+$. The argument for $\overline{\nu}_-$ is almost
verbally the same. Since $\overline{\nu}_+$ is isotropic, we can
find $n + 1$ linearly independent vectors in its support, say
$\{w_1,\dots,w_{n+1}\}$. If $w_0=\sum_{i=1}^{n+1} c_i w_i$ is an
arbitrary vector in $\mathrm{supp}\,\overline\nu_+$ such that,
without loss of generality, $c_1 \neq 0$, then by the equality
conditions of the Ball--Barthe Lemma,
\begin{align*}
\phi_{w_0}(w_0 \cdot z)&\phi_{w_2}(w_2 \cdot z)\cdots \phi_{w_{n+1}}(w_{n+1} \cdot z) \\
   &= \phi_{w_1}(w_1 \cdot z)\phi_{w_2}(w_2 \cdot z)\cdots \phi_{w_{n+1}}(w_{n+1} \cdot z)
\end{align*}
for every $z\in D$. Since $\phi_w$ is positive for every $w \in
\mathrm{supp}\,\overline\nu_+$, evaluating partial derivatives
with respect to $z$ at $\lambda e_{n+1}$ yields
\[\phi'_{w_0}(w_0 \cdot \lambda e_{n+1})w_0=\phi'_{w_1}(w_1 \cdot \lambda e_{n+1})w_1\]
for every $\lambda >0$. By the remark after Lemma
\ref{le:IsotrChar}, the support of $\overline{\nu}_+$ cannot
contain two antipodal points. Therefore, we either have that
$w_0=w_1$ or $\phi'_{w_0}(w_0 \cdot \lambda e_{n+1})=0 \text{ for
all }\lambda
>0.$ Since $w_0 \cdot e_{n+1} > 0$, the latter implies that
$\phi_{w_0}$ is constant, a contradiction. Consequently,
$\mathrm{supp}\,\overline{\nu}_+ = \{w_1,\dots,w_{n+1}\}$. Since
$\overline\nu$ is isotropic, it is easy to see (cf.\
\textbf{\cite{LYZ2004}}) that $w_1,\dots,w_{n+1}$ must be an
orthonormal basis of $\mathbb{R}^{n+1}$.

From the definition of $\overline{\nu}_+$, it follows that
$\mathrm{supp}\,\nu = \{u_1, \ldots, u_{n+1}\}$, where
$g_+(u_i)/\|g_+(u_i)\| = w_i$ for $1 \leq i \leq n + 1$. Using
definition (\ref{constemb}) of $g_+$, we obtain
\[0=w_i\cdot w_j=\frac{(u_i,\frac{1}{\sqrt n})\cdot (u_j,\frac{1}{\sqrt n})}{\sqrt{1+\frac{1}{n}}\sqrt{1+\frac{1}{n}}},\quad 1\leq i\neq j\leq n+1.\]
In other words, $u_i\cdot u_j=-\frac{1}{n}$ for all $i\neq j$.
Hence, $\mathrm{conv}\,\mathrm{supp}\,\nu$ must be a regular
simplex. \hfill $\blacksquare$

\vspace{1cm}

\centerline{\large{\bf{ \setcounter{abschnitt}{5}
\arabic{abschnitt}. Proof of the main results}}} \reseteqn
\alpheqn \setcounter{theorem}{0}

\vspace{0.6cm}

After these preparations, we are now in a position to give the
proofs of Theorems \ref{main1} and \ref{main2}. We start with the
following refinement of Theorem 1:

\begin{theorem} \label{main1const}
  Suppose $f$ is a positive continuous function on $S^{n-1}$ and $\nu$ is an isotropic $f$-centered measure.
  Then
\begin{equation} \label{mainineq1}
V(W_{\nu,f}) \leq \frac{ (n+1-\mathrm{disp}\,W_{\nu,f})^{n+1}}{n!(n+1)^{(n+1)/2}}\| f\|^{n}_{L^2(\nu)}
\end{equation}
with equality if and only if $\mathrm{conv}\,\mathrm{supp}\,\nu$
is a regular simplex inscribed in $S^{n-1}$ and $f$ is constant
on $\mathrm{supp}\,\nu$.
\end{theorem}

\noindent {\it Proof.} By the definition of $W_{\nu,f}$ and
$\mathrm{disp}\,W_{\nu,f}$, we may assume $\|f\|_{L^2(\nu)}=1$.
Let $\overline\nu$ denote the measure on $S^n$ defined by
(\ref{isoemb}), isotropically embedded by $g_-(u)=(-u,f(u))$, $u
\in S^{n-1}$ (here we use Lemma \ref{le:IsotrChar}).

\pagebreak

Next, let $C \subseteq \mathbb{R}^{n+1}$ denote the cone with apex
at the origin defined by
\[C=\bigcup_{r>0} r W_{\nu,f}\times \{r\} \subseteq \mathbb R^{n+1}.\]
Clearly, $e_{n+1} \in C$. Moreover, since $w \in
\mathrm{supp}\,\overline\nu \subseteq \mathbb R^{n}\times\mathbb
R$ if and only if
\begin{equation} \label{suppnudarst}
 w=\frac{(-u,f(u))}{\sqrt{1+f(u)^2}}
\end{equation}
 for some $u \in \mathrm{supp}\, \nu$, we have
that, for every $w \in\mathrm{supp}\,\overline\nu$ and $z =
(rx,r) \in C$,
\[w\cdot z=\frac{-u\cdot rx+rf(u)}{\sqrt{1+f(u)^2}} \geq 0.\]

For $w\in\mathrm{supp}\,\overline\nu$, define the smooth and
strictly increasing function \linebreak $T_w: (0,\infty)
\rightarrow \mathbb{R}$ by
\[\int_{-\infty}^{T_w(t)}e^{-\pi s^2}\, ds = \frac{1}{e_{n+1}\cdot w}\int_0^t \exp \left ( -\frac{s}{e_{n+1}\cdot w} \right ) \, ds.\]
Differentiating both sides with respect to $t$ yields
\[T_w'(t)\,e^{-\pi T_w(t)^2} = \frac{1}{e_{n+1}\cdot w}\,\exp \left ( -\frac{t}{e_{n+1}\cdot w} \right ).\]
Taking the $\log$ of both sides and putting $t=w\cdot z$ for
$w\in\mathrm{supp}\,\overline\nu$ and $z\in \mathrm{int}\, C$, we
obtain
\begin{equation}\label{eq:ModWulffMongeAmp}
\log T_w'(w \cdot z)-\pi\, T_w(w\cdot z)^2=-\log(e_{n+1}\cdot
w)-\frac{w\cdot z}{e_{n+1}\cdot w}.
\end{equation}

Furthermore define $T:\mathrm{int}\,C \rightarrow
\mathbb{R}^{n+1}$ by
\[T(z)=\int_{S^{n}} T_w(w\cdot z)\,w \, d \overline\nu(w).\]
A straightforward computation yields that, for every $z \in
\mathrm{int}\,C$,
\begin{equation} \label{dT}
dT(z)=\int_{S^{n}} T_w'(w\cdot z)\,w\otimes w\, d \overline\nu(w).
\end{equation}
Since $T_w'$ is a positive function, it follows that the matrix
$dT(z)$ is positive definite for $z\in \mathrm{int}\, C$.
Consequently, $T$ is a diffeomorphism onto its image. Moreover,
since
\[\|T(z)\|^2  = \int_{S^{n}} T_{w}(w\cdot z) (T(z)\cdot
w)\,d\overline\nu(w)\] and $\overline{\nu}$ is isotropic, we
obtain from an application of the Cauchy--Schwarz inequality that
\begin{align}\label{eq:CS}
\|T(z)\|^2 & \leq  \int_{S^{n}} T_{w}(w \cdot z)^2 \, d
\overline\nu(w).
  \end{align}

Now, by \eqref{eq:ModWulffMongeAmp}, followed by an application of
the Ball--Barthe Lemma with $t(w)=T_w'(w \cdot z)$, (\ref{dT}),
\eqref{eq:CS}, and a change of variables it follows that
  \begin{align*}
       \exp &  \left ( - \int_{S^{n}} \log(e_{n+1}\cdot w)\, d \overline\nu(w) \right ) \, \int_{\mathrm{int}\, C}\exp \left ( \int_{S^{n}} -\frac{w \cdot z}{e_{n+1}\cdot w} \, d \overline\nu(w) \right) \, d  z \\
          &=\int_{\mathrm{int}\, C} \exp \left ( \int_{S^{n}} \log T_w'(w \cdot z) \, d \overline\nu(w)\right )\exp\left(-\pi\int_{S^{n}} T_w(w \cdot z)^2\, d \overline\nu(w)\right)\, d  z \\
          &\leq \int_{\mathrm{int}\, C}\det \, dT(z)\exp\left(-\pi\|T(z)\|^2\right)\, d  z \leq \int_{\mathbb{R}^{n+1}} \exp\left(-\pi\|z\|^2\right)\, d
          z=1.
  \end{align*}
Equivalently, by the definition of the cone $C$,
    \begin{align}\begin{split}\label{eq:NormBLApp}
        \int_0^\infty \int_{r W_{\nu,f}} \exp\left(\int_{S^{n}}-\frac{(x,r)\cdot w}{e_{n+1}\cdot w}\, d  \overline\nu (w) \right)&\, d  x\, d  r \\
        \leq \exp \bigg (\int_{S^{n}} & \log \left ( (e_{n+1}\cdot w)^2 \right )\, d \overline\nu(w)\bigg )^{1/2}.
    \end{split}\end{align}

Since $\overline{\nu}$ is isotropic, an application of Jensen's
inequality to the right-hand side of \eqref{eq:NormBLApp} yields
\begin{equation} \label{zwischen1}
\int_0^\infty \int_{r W_{\nu,f}}
\exp\left(\int_{S^{n}}-\frac{(x,r)\cdot w}{e_{n+1}\cdot w}\, d
\overline\nu (w) \right) \, d  x\, d  r \leq
\left(\frac{1}{n+1}\right)^{(n+1)/2}.
\end{equation}
In order to obtain the desired inequality (\ref{mainineq1}), it
remains to show that the left-hand side of (\ref{zwischen1})
dominates
\[ V(W_{\nu,f})\, n!\, (n+1 - \mathrm{disp}\, W_{\nu,f})^{-(n+1)}.   \]

To see this, first note that, by definition (\ref{isoemb}) of
$\overline{\nu}$, the fact that we use the embedding
$g(u)=(-u,f(u))$ and since $\nu$ is $f$-centered, the left-hand
side of (\ref{zwischen1}) is equal to
\begin{equation} \label{zwischen2}
\int_0^\infty e^{-(n+1)r}\int_{r W_{\nu,f}}
\exp\left(\int_{S^{n-1}} \frac{x\cdot u}{f(u)}\, d
\nu(u)\right)\, dx \, dr.
\end{equation}

Since
\begin{equation} \label{dispdarst}
 r\, \mathrm{disp}\, W_{\nu,f} = \frac{1}{V(r W_{\nu,f})}\int_{r
W_{\nu,f}}\int_{S^{n-1}} \frac{x\cdot u}{f(u)}\, d \nu(u)\,dx,
\end{equation}
another application of Jensen's inequality, yields
\[\int_{r W_{\nu,f}}
\exp\left(\int_{S^{n-1}} \frac{x\cdot u}{f(u)}\, d
\nu(u)\right)\, dx \geq V(rW_{\nu,f})\,e^{r\, \mathrm{disp}\,
W_{\nu,f}}.   \] Consequently, (\ref{zwischen2}) is larger than
\[V(W_{\nu,f}) \int_0^\infty r^n e^{-(n+1 - \mathrm{disp}\, W_{\nu,f})r}\, dr = V(W_{\nu,f})\, n!\, (n+1 - \mathrm{disp}\, W_{\nu,f})^{-(n+1)},  \]
where we have used that $\mathrm{disp}\, W_{\nu,f} \leq n$, which
follows easily from (\ref{dispdarst}) and the definition of
$W_{\nu,f}$. This completes the proof of inequality
(\ref{mainineq1}).

Assume now that there is equality in inequality
(\ref{mainineq1}). By the equality conditions of Jensen's
inequality, equality in (\ref{zwischen1}) can hold only if
$e_{n+1} \cdot w$ is constant for every $w \in
\mathrm{supp}\,\overline\nu$. Since any $w \in
\mathrm{supp}\,\overline\nu$ is of the form (\ref{suppnudarst}),
this implies that $f$ is constant on the support of $\nu$. By the
normalization $\|f\|_{L^2(\nu)}=1$, we must have
$f\equiv\frac{1}{\sqrt n}$ on $\mathrm{supp}\,\nu$. Now it is
easy to check that the assumptions of Lemma \ref{equcond}, where
$D = \mathrm{int}\,C$ and $\phi_w=T'_w$, are satisfied. Hence, an
application of Lemma \ref{equcond} concludes the proof. \hfill
$\blacksquare$

\vspace{0.3cm}

In order to establish Theorem \ref{main2}, we will use a transport
map $\widehat{T}_w$ that is in some sense dual to the function
$T_w$ used in the proof above.

\vspace{0.3cm}

\noindent {\it Proof of Theorem \ref{main2}.} As before we may
assume that $\|f\|_{L^2(\nu)}=1$. In the following we denote by
$\overline\nu$ the measure $S^n$ defined by (\ref{isoemb}),
isotropically \linebreak embedded by $g_+(u)=(u,f(u))$, $u \in
S^{n-1}$ (where we again use Lemma \ref{le:IsotrChar}).

For $w\in\mathrm{supp}\,\overline\nu$, define the smooth and
strictly increasing function \linebreak $\widehat{T}_w: \mathbb{R}
\rightarrow (0,\infty)$ by
\[e_{n+1}\cdot w \int_0^{\widehat{T}_w(t)} e^{-s(e_{n+1}\cdot w)} \, d  s =\int_{-\infty}^t e^{-\pi s^2}\, d  s.\]
Differentiating both sides with respect to $t$ yields
\[(e_{n+1}\cdot w)\,\widehat{T}_w'(t)\,e^{-\widehat{T}_w(t)(e_{n+1}\cdot w)} =  e^{-\pi t^2}.  \]
Taking the $\log$ of both sides and putting $t=w\cdot z$ for
$w\in\mathrm{supp}\,\overline\nu$ and $z\in \mathbb{R}^{n+1}$, we
obtain
  \begin{equation}\label{eq:PolarWulffMongeAmp}
       \log \widehat{T}_w'(w \cdot z)=\widehat{T}_w(w \cdot z)(e_{n+1}\cdot w) -\pi (w \cdot z)^2 - \log(e_{n+1}\cdot w).
  \end{equation}
Define the map $\widehat{T}:\mathbb R^{n+1}\to\mathbb R^{n+1}$ by
  \[\widehat{T}(z):=\int_{S^{n}} \widehat{T}_w(w\cdot z)\, w\, d \overline\nu(w).  \]
 The Jacobian of $\widehat T$ is given by
  \[d\widehat T(z)=\int_{S^{n}} \widehat{T}_w'(w \cdot z)\, w\otimes w\, d \overline\nu(w),  \]
  which shows that $d\widehat T(z)$ is a positive definite matrix for every $z \in \mathbb{R}^{n+1}$. Consequently, $\widehat T$ is a diffeomorphism onto its image.
  In fact, we claim that its image is contained in the cone
  \[\widehat C:=\bigcup_{r>0} r W_{\nu,f}^*\times\{r\}.\]
  To prove this, we have to show that if $\widehat T(z) = (x,r) \in \mathbb{R}^{n+1}=\mathbb{R}^n \times \mathbb{R}$ and $y \in
  W_{\nu,f}$, then $x \cdot y \leq r$. By definition (\ref{isoemb}) of $\overline\nu$, the definition
  of $\widehat{T}$ and the fact that $u\cdot y \leq f(u)$ for every $u \in \mathrm{supp}\,\nu$, we obtain
  \begin{align*}
      x \cdot y & = \int_{S^{n-1}} \widehat{T}_w\left(\frac{(u,f(u))}{\sqrt{1+f(u)^2}}\cdot z\right) (u \cdot y)\sqrt{1+f(u)^2}\, d \nu(u)\\
           &\leq \int_{S^{n-1}} \widehat{T}_w\left(\frac{(u,f(u))}{\sqrt{1+f(u)^2}} \cdot z\right) f(u)\sqrt{1+f(u)^2}\, d \nu(u) \\
           &=\int_{S^{n}} \widehat{T}_w(w\cdot z) (e_{n+1}\cdot w) \, d \overline\nu(w)=r.
  \end{align*}
Since
\[n!\,V(W_{\nu,f}^*)= \int_0^\infty \int_{r W_{\nu,f}^*} e^{-e_{n+1}\cdot z} \, d  x\, d  r = \int_{\widehat{C}} e^{-e_{n+1}\cdot z} \, d  z, \]
a change of variables, followed by an application of the
Ball--Barthe Lemma, \eqref{eq:PolarWulffMongeAmp} and the fact
that $\overline \nu$ is isotropic, yields
  \begin{align*}\begin{split}
n!\,V(W_{\nu,f}^*) & \geq \int_{\mathbb{R}^{n+1}} e^{-e_{n+1}\cdot \widehat T(z)}\det\, d \widehat T(z)\, d  z\\
      &\geq \int_{\mathbb{R}^{n+1}} e^{-e_{n+1}\cdot \widehat T(z)}\exp\left(\int_{S^{n}} \log \widehat T'_w(w \cdot z)\, d \overline\nu(w)\right)\, d  z\\
      &= \exp\left(-\int_{S^{n}} \log(e_{n+1}\cdot w)\, d  \overline\nu(w)\right)\int_{\mathbb{R}^{n+1}} \exp\left(-\pi\|z\|^2\right)\, dz \\
      & = \exp\left(\frac{1}{n+1}\int_{S^{n}} \log\left ( (e_{n+1}\cdot w)^2 \right )\, d  \overline\nu(w)\right)^{-(n+1)/2}.\end{split}
  \end{align*}
Thus, using Jensen's inequality and again the fact that
$\overline \nu$ is isotropic, we obtain the desired inequality
\begin{align} \label{mainineq2}
n!\,V(W_{\nu,f}^*) \geq (n+1)^{(n+1)/2}.
\end{align}

Assume now that there is equality in inequality (\ref{mainineq2}). As in the proof of Theorem \ref{main1const} we conclude that this is possible only if $f$ is constant on the support of $\nu$. Thus, by the normalization $\|f\|_{L^2(\nu)}=1$, we must have $f \equiv \frac{1}{\sqrt{n}}$ on $\mathrm{supp}\,\nu$. In order to apply Lemma \ref{equcond}, define the open cone $D \subseteq \mathbb{R}^{n+1}$ by
\[D=\{z\in\mathbb R^{n+1}:  w \cdot z > 0 \text{ for every } w\in\mathrm{supp}\, \overline\nu\}.\]
Clearly, the assumptions of Lemma \ref{equcond}, where $\phi_w = \widehat{T}_w'$ are now satisfied. Hence, an application of Lemma \ref{equcond} concludes the proof. \hfill $\blacksquare$

\vspace{1cm}

\centerline{\large{\bf{ \setcounter{abschnitt}{6}
\arabic{abschnitt}. Applications}}} \reseteqn \alpheqn
\setcounter{theorem}{0}

\vspace{0.6cm}

In this final section, we show how Theorem 1 and 2 directly imply \linebreak previously established reverse isoperimetric inequalities which have simplices as extremals.
We begin with consequences of Theorem 1, such as Ball's \linebreak volume ratio inequality and its $L_2$ analog by Lutwak, Yang, and Zhang, and conclude this section with
dual results (including Corollary \ref{introgam2}), which can be deduced from Theorem 2.

First suppose that $\nu$ is an 1-centered isotropic measure on
$S^{n-1}$. Then
\[ W_{\nu,1}=(\mathrm{conv}\,\mathrm{supp}\,\nu)^*\]
and $\mathrm{disp}\, W_{\nu,1}=0$. Consequently, Theorem \ref{main1} reduces to the following result of Lutwak, Yang and Zhang \textbf{\cite{LYZ2007}}.

\pagebreak

\begin{corol}\label{cor:PolarConvSupp} {\bf (\!\!\cite{LYZ2007})}
  If $\nu$ is an 1-centered isotropic measure on $S^{n-1}$, then
\begin{equation} \label{cor1}
V((\mathrm{conv}\,\mathrm{supp}\,\nu)^*)\leq \frac{n^{n/2}(n+1)^{(n+1)/2}}{n!},
\end{equation}
with equality if and only if $\mathrm{conv}\,\mathrm{supp}\,\nu$ is a regular simplex inscribed in $S^{n-1}$.
\end{corol}

Ball \textbf{\cite{ball91b}} had first established inequality (\ref{cor1}), but without the equality conditions.
For discrete measures, these were obtained by Barthe
\textbf{\cite{barthe98i}}.

Corollary \ref{cor:PolarConvSupp} allows for a geometric
interpretation, known as Ball's \linebreak volume ratio
inequality, which gives an upper bound for the ratio between the
volume of a convex body and its John ellipsoid:

\vspace{0.4cm}

\noindent {\bf Ball's Volume Ratio Inequality (\!\!\cite{ball91b, barthe98i})} \emph{If $K \subseteq \mathbb{R}^n$ is a convex body, then
\begin{equation}\label{eq:VRIneq}
 \frac{V(K)}{V(\mathrm{J}K)}\leq
 \frac{n^{n/2}(n+1)^{(n+1)/2}}{\kappa_nn!},
\end{equation}
with equality if and only if $K$ is a simplex.}

\vspace{0.2cm}

\noindent {\it Proof.} Without loss of generality we may assume
that $K$ is in John position, \linebreak that is $\mathrm{J}K=B$.
By Proposition \ref{johnth}, there exists an 1-centered isotropic
measure $\nu$ on $S^{n-1}$ supported by contact points of $K$ and
$B$. Clearly, $K\subseteq
W_{\nu,1}=(\mathrm{conv}\,\mathrm{supp}\,\nu)^*$. Thus, Corollary
\ref{cor:PolarConvSupp} implies (\ref{eq:VRIneq}) along with its
equality conditions. \hfill $\blacksquare$

\vspace*{0.3cm}

If  $K \subseteq \mathbb{R}^n$ is a convex body such that
$\mathrm{J}K$ is centered at the origin, then we can replace
$\mathrm{J}K$ in inequality \eqref{eq:VRIneq} by the $L_\infty$
John ellipsoid $\mathrm{E}_\infty K$. A combination of
(\ref{eq:VRIneq}) and \eqref{eq:VIneqpq} thus yields the
following {\it $L_p$ volume ratio inequality} for the entire
family of $L_p$ John ellipsoids $\mathrm{E}_pK$,
$0<p\leq\infty$\,:
      \begin{equation}\label{eq:VRLp}
          \frac{V(K)}{V(\mathrm{E}_p K)}\leq \frac{n^{n/2}(n+1)^{(n+1)/2}}{\kappa_nn!},
      \end{equation}
with equality if and only if $K$ is a simplex with centroid at
the origin.

Using a different approach, Lutwak, Yang and Zhang established
the case $p=2$ of inequality \eqref{eq:VRLp} in
\textbf{\cite{LYZ-newellassconbod}}. This $L_2$ volume ratio
inequality is also a {\it direct} consequence of Theorem
\ref{main1}, where in addition we can replace the assumption that
$\mathrm{J}K$ is centered at the origin by the more natural
assumption that $K$ has centroid at the origin:

\begin{corol} If $K \subseteq \mathbb{R}^n$ is a convex body with centroid at the origin, then
  \begin{equation}\label{eq:L2VRIneq}
              \frac{V(K)}{V(\mathrm{E}_2 K)}\leq \frac{n^{n/2}(n+1)^{(n+1)/2}}{\kappa_nn!},
          \end{equation}
with equality if and only if $K$ is a simplex with centroid at
the origin.
\end{corol}
\noindent {\it Proof.}
 Without loss of generality we may assume that $K$ is in $L_2$ John position, that is $\mathrm{E}_2K=B$.
By Proposition \ref{lpjohnth}, this implies that the measure $\nu:=\frac{1}{V(K)}S_2(K,\cdot)$
    is isotropic. Moreover, by (\ref{eq:SurfCentred}) and (\ref{absc}), the measure $\nu$ is $h(K,\cdot)$-centered.
    By \eqref{eq:KIntersectionHalfSp}, \eqref{vol1} and \eqref{defvp}, we have
    \[ W_{\nu,h(K,\cdot)}=K\qquad\text{ and }\qquad \|h(K,\cdot)\|_{L^2(\nu)}=\sqrt n. \]
Thus, $\mathrm{disp}\, W_{\nu,h(K,\cdot)}=0$ and Theorem \ref{main1} implies inequality (\ref{eq:L2VRIneq}) along with its equality conditions. \hfill $\blacksquare$

\vspace{0.3cm}

A combination of (\ref{eq:L2VRIneq}) and (\ref{eq:VIneqpq}) shows
that inequality \eqref{eq:L2VRIneq} (under the assumption that
$\mathrm{cd}\,K = o$) holds true if $\mathrm{E}_2K$ is replaced by
$\mathrm{E}_pK$, $0<p\leq 2$.

\vspace{0.2cm}

We now turn to special cases of Theorem 2. To this end, it is
useful to keep the following (easily verified) alternative
representation of the polar Wulff shape $W_{\nu,f}^*$ of a given
$f$-centered isotropic measure $\nu$ in mind:
\begin{equation} \label{wulffpolar}
W_{\nu,f}^*=\mathrm{conv}\, \left\{\frac{u}{f(u)}:\;
u\in\mathrm{supp}\,\nu\right\}.
\end{equation}

If $\nu$ is now an $1$-centered isotropic measure on $S^{n-1}$,
then, by (\ref{wulffpolar}), \linebreak Theorem \ref{main2}
reduces to the following result of Lutwak, Yang, Zhang
\textbf{\cite{LYZ2007}}.

\begin{corol}\label{cor:ConvSupp} {\bf (\!\!\cite{LYZ2007})}
  If $\nu$ is an $1$-centered isotropic measure on $S^{n-1}$, then
\begin{equation} \label{dualnuVR}
V(\mathrm{conv}\,\mathrm{supp}\,\nu)\geq
\frac{(n+1)^{(n+1)/2}}{n^{n/2}n!},
\end{equation}
with equality if and only if $\mathrm{conv}\,\mathrm{supp}\,\nu$
is a regular simplex inscribed in $S^{n-1}$.
\end{corol}

For discrete measures, Corollary \ref{cor:ConvSupp} was first
established by Barthe \textbf{\cite{barthephd}}. A more geometric
reformulation of inequality (\ref{dualnuVR}) is a dual result to
Ball's volume ratio inequality:

\pagebreak

\noindent {\bf Barthe's Dual Volume Ratio Inequality
(\!\!\cite{barthephd})} \emph{If $K \subseteq \mathbb{R}^n$ is a
convex body, then
  \begin{equation}\label{eq:DualVRIneq}
    V(K^*)V(\mathrm{J}K)\geq \frac{(n+1)^{(n+1)/2}\kappa_n}{n^{n/2}n!},
  \end{equation}
with equality if and only if $K$ is a simplex with centroid at
the origin.}

\vspace{0.2cm}

\noindent {\it Proof.} First, we use Corollary \ref{cor:ConvSupp}
to deduce Barthe's \textbf{\cite{barthephd}} outer volume ratio
inequality
\begin{equation}\label{eq:outerVRIneq}
\frac{V(K)}{V(\mathrm{L}K)}\geq
\frac{(n+1)^{(n+1)/2}}{n^{n/2}n!\kappa_n},
\end{equation}
where $\mathrm{L}K$ denotes the Loewner ellipsoid of $K$, that is
the ellipsoid of minimal volume containing $K$. To this end, we
may assume without loss of generality that $\mathrm{L}K=B$. Then,
there exists an $1$-centered isotropic measure $\nu$ on $S^{n-1}$
supported by contact points of $K$ and $B$ (compare the remark
after \linebreak Proposition \ref{johnth}). Clearly,
$\mathrm{conv}\,\mathrm{supp}\,\nu\subseteq K$. Thus, using
Corollary \ref{cor:ConvSupp}, we obtain inequality
(\ref{eq:outerVRIneq}) with equality if and only if $K$ is a
simplex.

Using the definitions of $\mathrm{J}K$ and $\mathrm{L}K$ and the
fact that for every ellipsoid $E$ containing the origin in its
interior, $V(E)V(E^*)\geq \kappa_n^2$ with equality precisely for
origin-symmetric ellipsoids, we obtain
\[V(\mathrm{J}K)V(\mathrm{L}K^*) \geq V((\mathrm{L}K^*)^*)V(\mathrm{L}K^*) \geq \kappa_n^2. \]
Combining this with inequality \eqref{eq:outerVRIneq}, where $K$
is replaced by $K^*$, yields the desired inequality
\eqref{eq:DualVRIneq} along with its equality conditions. \hfill
$\blacksquare$

\vspace{0.3cm}

If $K \subseteq \mathbb{R}^n$ is a convex body such that
$\mathrm{J}K$ is centered at the origin, that is $\mathrm{J}K =
\mathrm{E}_\infty K$, then a combination of (\ref{eq:DualVRIneq})
and \eqref{eq:VIneqpq} yields the following dual to inequality
(\ref{eq:VRLp}) for $0<p\leq\infty$\,:
\begin{equation}\label{eq:DualVRIneqLp}
  V(K^*)V(\mathrm{E}_pK)\geq \frac{(n+1)^{(n+1)/2}\kappa_n}{n^{n/2}n!},
\end{equation}
with equality if and only if $K$ is a simplex with centroid at
the origin.

It is an important task in convex geometry to find sharp lower
bounds \linebreak for the volume of $K^*$ in terms of other
natural geometric quantities of $K$ \linebreak (see e.g.\
\textbf{\cite{barthefrad11, BM-newvolratproconsymbodbfrn,
K-mahcontogaulinint, NPRZ-remmahconlocminunicub, R-zonminvol}}
for results in this direction). Unfortunately, for inequality
(\ref{eq:DualVRIneqLp}), the requirement that the John ellipsoid
of $K$ is centered at the origin cannot, in general, be omitted
for $p > 2$. For the $L_2$ case, Lutwak, Yang, and Zhang
\textbf{\cite{LYZ2009}} have discovered that this additional
assumption is in fact unnecessary. Their result was stated in the
introduction as Corollary \ref{introgam2} and is (in the
following equivalent formulation) also a direct consequence of
Theorem 2:

\begin{corol} {\bf (\!\!\cite{LYZ2009})}
  If $K \subseteq \mathbb{R}^n$ is a convex body, then
    \begin{equation}\label{eq:DualVRIneqL2}
    V(K^*)V(\mathrm{E}_2K)\geq \frac{(n+1)^{(n+1)/2}\kappa_n}{n^{n/2}n!},
    \end{equation}
with equality if and only if $K$ is a simplex with centroid at
the origin.
\end{corol}

\noindent {\it Proof.} Without loss of generality we may assume
that $K$ is in $L_2$ John position. By Proposition
\ref{lpjohnth}, this implies that the $h(K,\cdot)$-centered
measure $\nu:=\frac{1}{V(K)}S_2(K,\cdot)$ is isotropic. By
\eqref{eq:KIntersectionHalfSp}, \eqref{vol1} and \eqref{defvp},
we have
    \[ W_{\nu,h(K,\cdot)}=K\qquad\text{ and }\qquad \|h(K,\cdot)\|_{L^2(\nu)}=\sqrt n. \]
Thus, Theorem \ref{main2} reduces to the desired statement.
\hfill $\blacksquare$

\vspace{0.3cm}

A combination of (\ref{eq:DualVRIneqL2}) and (\ref{eq:VIneqpq})
again shows that inequality \eqref{eq:DualVRIneqL2} remains true
if $\mathrm{E}_2K$ is replaced by $\mathrm{E}_pK$, $0<p\leq 2$.

\vspace{0.1cm}

We finally remark that all the special cases of Theorem 1 and 2
presented in this section have natural analogues for
origin-symmetric convex bodies, where the cube instead of the
simplex plays the extremal role. These are of course special
cases of the volume inequalities for symmetric Wulff shapes
stated in Section 3.

\vspace{0.5cm}

\noindent {{\bf Acknowledgments} The work of the authors was
supported by the \linebreak Austrian  Science Fund (FWF), within
the project ``Minkowski valuations and geometric inequalities",
Project number: P\,22388-N13.

\begin{small}

Vienna University of Technology \par Institute of Discrete
Mathematics and Geometry \par Wiedner Hauptstra\ss e 8--10/1046
\par A--1040 Vienna, Austria

franz.schuster@tuwien.ac.at \par

manuel.weberndorfer@tuwien.ac.at

\end{small}

\end{document}